\documentclass[12pt]{amsart}
\usepackage[mathscr]{eucal}
\usepackage{amssymb}
\usepackage{latexsym}
\usepackage{amsthm}
\usepackage{amsmath}
\usepackage{graphicx}

\theoremstyle{plain}

\newtheorem{thm}{Theorem}[section]
\newtheorem{cor}{Corollary}[section]
\newtheorem{lem}{Lemma}[section]

\theoremstyle{definition}

\newtheorem{prf}{Proof}[section]

\DeclareMathOperator*{\esssup}{ess\, sup}

\title{Functional analytic approach to Ces\`{a}ro mean}
\author{Ryoichi Kunisada}

\address{Faculty of Education and Integrated Arts and Science, Waseda University, Shinjuku-ku, Tokyo 169-8050, Japan}

\email{rkunisada@aoni.waseda.jp}

\date{}

\begin{document}
\maketitle

\begin{abstract}
We study a certain class $\mathcal{P}$ of positive linear functionals $\varphi$ on $L^{\infty}([1,\infty))$ for which $\varphi(f) = \alpha$ 
if $\lim_{x \to \infty} \frac{1}{x} \int_1^x f(t)dt = \alpha$. It turns out that translations $f(x) \mapsto f(rx)$ on $L^{\infty}([1, \infty))$, where $r \in [1, \infty)$, which are induced by the action of the multiplicative semigroup $[1, \infty)$ on itself, plays an intrinsic role in the study of $\mathcal{P}$. We also deal with an analogue $\mathcal{K}$ of $\mathcal{P}$ of positive linear functionals on $L^{\infty}([0, \infty))$ partaining to the action of
the additive semigroup $[0, \infty)$ on itself. In particular, we give some expressions of maximal possible values of $\mathcal{P}$ and $\mathcal{K}$ for a given function respectively.
\end{abstract}

\bigskip

\section{Introduction}
Let $L^{\infty}(\mathbb{R}_+^{\times})$ be the space of all essentially bounded functions on $\mathbb{R}_+^{\times} = [1, \infty)$. Let $L^{\infty}(\mathbb{R}_+^{\times})^*$ be the dual space of $L^{\infty}(\mathbb{R}_+^{\times})$. Given $f(x) \in L^{\infty}(\mathbb{R}_+^{\times})$, define
\[M(f) = \lim_x \frac{1}{x} \int_1^x f(t)dt \]
if this limit exists. Obviously, it is an integral analogy of the notion of Ces\`{a}ro mean.
Let $E$ be the space of all functions in $L^{\infty}(\mathbb{R}_+^{\times})$ having the limit $M(f)$ and $E_0$ be its subspace consisting of all $f \in E$ for which $M(f) = 0$.

Our main interest of this paper is the set $\mathcal{C}$ of normalized positive linear functionals on $L^{\infty}(\mathbb{R}_+^{\times})$ which vanish on $E_0$, that is, $\mathcal{C} = E_0^{\bot} \cap S^+_{L^{\infty}(\mathbb{R}_+^{\times})^*}$, where $E_0^{\bot}$ is the annihilator of $E_0$ and $S^+_{L^{\infty}(\mathbb{R}_+^{\times})^*}$ is the positive part of the unit sphere of $L^{\infty}(\mathbb{R}_+^{\times})^*$. In other words, $\varphi \in \L^{\infty}(\mathbb{R}_+^{\times})^*$ is an element of $\mathcal{C}$ if and only if it is an extension of the functional $M : E \rightarrow \mathbb{R}$. 

Of particular importance for our study is the following sublinear functional
\begin{align}
\overline{P}(f) &= \lim_{\theta \to 1^+} \limsup_{x \to \infty} \frac{1}{\theta x - x} \int_x^{\theta x} f(t)dt \notag \\
&= \lim_{\theta \to 1^+} \limsup_{x \to \infty} \frac{1}{\theta - 1} \int_1^{\theta} f(xt)dt,  \quad f(x) \in L^{\infty}(\mathbb{R}_+^{\times}), \notag
\end{align}
which turns out to give the maximal value attained by elements of $\mathcal{C}$ for each $f \in L^{\infty}(\mathbb{R}_+^{\times})$. That is, 
\[\sup_{\varphi \in \mathcal{C}} \varphi(f) = \overline{P}(f) \]
holds. One of our main aim of this paper is to prove this assertion by functional analytic methods and to give another expression of $\sup_{\varphi \in \mathcal{C}} \varphi(f)$ which has a more simple form.

It should be noted that  these definitions derive from their summation counterparts, with which one is more familier. We give below a brief description of this notion. Let $l_{\infty}$ be the set of all real-valued bounded functions on natural numbers $\mathbb{N}$. Recall that for any $f \in l_{\infty}$, its Ces\`{a}ro mean $M_d(f)$ is define as
\[M_d(f) = \lim_{n \to \infty} \frac{1}{n} \sum_{i=1}^n f(i) \]
if this limit exists. Let $\mathcal{C}_d$ be the set of normalized positive linear functionals on $l_{\infty}$ which extend  Ces\`{a}ro mean. Also we define a sublinear functional $\overline{P}_d$ on $l_{\infty}$ as follows.
\[\overline{P}_d(f) = \lim_{\theta \to 1-} \limsup_{n \to \infty} \frac{\sum_{i \in [\theta n, n]} f(i)} {n - \theta n}, \]
where $f \in l_{\infty}$. Then it holds that
\[\sup_{\varphi \in \mathcal{C}_d} \varphi(f) = \overline{P}_d(f). \]
Based on the results of P\'{o}lya in [5], this is proved in [3]. Remark that elements of $\mathcal{C}_d$ and the functional $\overline{P}_d$ are sometimes called density measures and P\'{o}lya density respectively when they are ristricted to the characteristic functions of subsets of $\mathbb{N}$. Density measures have been studied by several authors, for example [1], [4], [6]. In particular, our main result mentioned above can be considered as an integral analogy of this result, which in fact, as we will see in Section 5, we can  deduce from the integral version with ease.

To this end, we introduce another class of continuous linear functionals on $L^{\infty}(\mathbb{R}_+)$ of the space of all essentially bounded functions on $\mathbb{R}_+ = [0, \infty)$. Given $f(x) \in L^{\infty}(\mathbb{R}_+)$, define
\[R(f) = \lim_{x \to \infty} \frac{1}{e^x} \int_0^x f(t)e^tdt \]
if this limit exists. Similarly as preceding paragraphs, let $F$ be the space of all functions $f$ in $L^{\infty}(\mathbb{R}_+)$ with the limit $R(f)$ exists and $F_0$ be its subspace consisting of all $f \in F$ for which $R(f) = 0$. Then we consider the set $\mathcal{D}$ of normalized positive functionals on $L^{\infty}(\mathbb{R}_+)$ which vanish on $F_0$, i.e., $\mathcal{D} = F_0^{\bot} \cap S^+_{L^{\infty}(\mathbb{R}_+)^*}$. It may be said that the idea of introducing the class $\mathcal{D}$ in connection with $\mathcal{C}$ is very natural. We can formulate a similar question of giving an expression of the maximal value of $\mathcal{D}$ for any fixed $f \in L^{\infty}(\mathbb{R}_+)$. Let us define a sublinear functional $\overline{K}$ on $L^{\infty}(\mathbb{R}_+)$ by
\begin{align} 
\overline{K}(f) &= \lim_{\theta \to 0^+} \limsup_{x \to \infty} \frac{1}{\theta} \int_x^{x+\theta} f(t)dt \notag \\
&= \lim_{\theta \to 0^+} \limsup_{x \to \infty} \frac{1}{\theta} \int_0^{\theta} f(t+x)dt, \quad f \in L^{\infty}(\mathbb{R}_+). \notag
\end{align}
Notice that the definition of $\overline{K}$ is analogous to that of $\overline{P}$ in the sense that we consider the action of the additive group of $\mathbb{R}$ in place of the action of the
multiplicative group $\mathbb{R}^{\times} = (0, \infty)$ of $\mathbb{R}$ in the definition of $\overline{P}$. Also, the additive counterpart of the preceding assertion about maximal values of $\mathcal{C}$ is given as follows:
\[\sup_{\varphi \in \mathcal{D}} \varphi(f) = \overline{K}(f) \]
holds for each $f \in L^{\infty}(\mathbb{R}_+)$. 

The paper is organaized as follows. After Section 2, which contains necessary notation and notions needed in the rest of the paper, we first consider the additive version in Section 3. We deal with the multiplicative version in Section 4, where arguments are similar to those in Section 3. Also we consider the relationship between $\mathcal{C}$ and $\mathcal{D}$ and show that there exists a natural affine homeomorphism between $\mathcal{C}$ and $\mathcal{D}$. Last section deal with the descrete version of resutls in Section 4. 

\section{Preliminaries}
Let $X$ be a set and $Y$ be a compact space and $f : X \rightarrow Y$ be a mapping. Let $\mathcal{U}$ be an ultrafilter on $X$. Then there exists an element $y$ of $Y$ such that $f^{-1}(U) \in \mathcal{U}$ holds for every neighborhood $U$ of $y$. In this case, we write $\mathcal{U}\mathchar`-\lim_x f(x) = y$ and say that $y$ is the limit of $f$ along $\mathcal{U}$.

Let $\mathbb{N}_0$ be the set of nonnegative integers and $\beta\mathbb{N}_0$ be its Stone-\v{C}ech compactification. Recall that $\beta\mathbb{N}_0$ is a compactification of $\mathbb{N}_0$ such that for any continuous mapping $\iota : \mathbb{N}_0 \rightarrow X$ of $\mathbb{N}_0$ into a compact space $X$, there exists a continuous extension $\overline{\iota} : \beta\mathbb{N}_0 \rightarrow X$ of $\iota$ to $\beta\mathbb{N}_0$.

We denote by $\tau : \beta\mathbb{N}_0 \rightarrow \beta\mathbb{N}_0$ the continuous extension of the mapping $\mathbb{N}_0 \ni n \mapsto n+1 \in \beta\mathbb{N}_0$. Let us denote $\mathbb{N}^* = \beta\mathbb{N}_0 \setminus \mathbb{N}_0$ and since the restriction of $\tau$ to $\mathbb{N}^*$ is a homeomorphism of $\mathbb{N}^*$ onto itself, the pair $(\mathbb{N}^*, \tau)$ is a topological dynamical system. 

Next we consider the compact space $\Omega$ of all maximal ideals of $C_{ub}(\mathbb{R}_+)$ of the space of all real-valued uniformly continuous bounded functions on $\mathbb{R}_+$, which can be viewed in a sense as a continuous version of $\beta\mathbb{N}_0$. It is known that the maximal ideal space $\Omega$ of $C_{ub}(\mathbb{R}_+)$ is given as follows (see [2] for further details of the following assertions): $\Omega = (\beta\mathbb{N}_0 \times [0 ,1]) / \sim$, the quotient space of the product space $\beta\mathbb{N}_0 \times [0,1]$, where $\sim$ is an equivalence relation on $\beta\mathbb{N}_0 \times [0 ,1]$ such that $\omega \sim \omega^{\prime}$ if and only if $ \omega = (\eta, 1)$ and $\omega^{\prime} = (\tau\eta, 0)$ for some $\eta \in \beta\mathbb{N}_0$. Then there exists an algebraic isomorphism $C_{ub}(\mathbb{R}_+) \ni f(x) \rightarrow \overline{f}(\omega) \in C(\Omega)$ of $C_{ub}(\mathbb{R}_+)$ onto $C(\Omega)$ of the set of all continuous functions on $\Omega$. Notice that for each $\omega = (\eta, t) \in \Omega$, the class of sets $\{t+A : A \in \eta \}$ is an ultrafilter on $\mathbb{R}_+$. In this way, we identify an element $\omega$ of $\Omega$ with an ultrafilter on $\mathbb{R}_+$ from now on. Then for every $f \in C_{ub}(\mathbb{R}_+)$, it holds that $\overline{f}(\omega) = \omega\mathchar`-\lim_x f(x)$.

Further, for each $s \ge 0$, let us define a continuous mapping $\tau^s : \Omega \rightarrow \Omega$ as follows:
\[\tau^s\omega = \tau^s(\eta, t) = (\tau^{[t+s]}\eta, t+s-[t+s]), \quad   \]
where $s \in \mathbb{R}$ and $[x]$ denotes the largest integer not exceeding a real number $x$. Let us denote $\Omega^* = \Omega \setminus \mathbb{R}$ and remark that the restriction of each $\tau^s$ to $\Omega^*$ is a homeomorphism of $\Omega^*$ onto itself. Thus the pair $(\Omega^*, \{\tau^s\}_{s \in \mathbb{R}})$ is a continuous flow.

\section{Extremal values of $\mathcal{D}$}
We consider the set of continuous linear functionals $\varphi$ on $L^{\infty}(\mathbb{R}_+)$ for which 
\[\varphi(f) \le \overline{K}(f) \]
holds for every $f \in L^{\infty}(\mathbb{R}_+)$ and denote it by $\mathcal{K}$.

Given $f(x) \in L^{\infty}(\mathbb{R}_+)$, we define a sequence $\{\tilde{f}_n(x)\}_{n \ge 0}$ of $L^{\infty}(\mathbb{T})$, where $\mathbb{T} = [0, 1]$, by
\[\tilde{f}_n(x) = f(x+n), \quad x \in [0, 1], \ n=0,1,2, \ldots. \]
Notice that $\{\tilde{f}_n(x)\}_{n \ge 0}$ is a uniformly bounded sequence of $L^{\infty}(\mathbb{T})$ and then it is a weak* relatively compact subset of $L^{\infty}(\mathbb{T})$. Then the mapping $\mathbb{N}^0 \ni n \mapsto \tilde{f}_n(x) \in L^{\infty}(\mathbb{T})$ can be extended continuously to $\beta\mathbb{N}^0$. For each $\eta \in \beta\mathbb{N}^0$, let us denote its image under this extended mapping by $\tilde{f}_{\eta}(x) \in L^{\infty}(\mathbb{T})$.
Notice that it can be expressed as the limit along an ultrafilter $\eta$, i.e., $\tilde{f}_{\eta}(x) = \eta\mathchar`-\lim_n \tilde{f}_n(x)$.

We need the lemma below which asserts that we can interchange limit and integral. 
Let $\{f_n(x)\}_{n \ge 0}$ be a uniformly bounded sequence of $L^{\infty}(\mathbb{T})$ and let
\[\eta\mathchar`-\lim_n f_n(x) = f_{\eta}(x) \]
for $\eta \in \mathbb{N}^* = \beta\mathbb{N}\setminus\mathbb{N}$, where the limit is taken with respect to the weak*-topology of $L^{\infty}(\mathbb{T})$. Now we consider their indefinite integrals; let us define
\[F_n(x) = \int_0^x f_n(t)dt, \quad x \in [0, 1], \ n=0, 1, 2, \ldots,  \]
\[G_{\eta}(x) = \int_0^x f_{\eta}(t)dt, \ x \in [0, 1]. \]
Then it is obviously that $\{F_n(x)\}_{n \ge 0}$ is a uniformly bounded and equicontinuous sequence of $C(\mathbb{T})$, the set of all continuous functions on $\mathbb{T}$. 
It means that $\{F_n(x)\}_{n \ge 0}$ is a relatively compact subset of $C(\mathbb{T})$ in its uniform topology. Therefore
we can consider for each $\eta \in \mathbb{N}^*$ the limit $F_{\eta}(x) = \eta\mathchar`-\lim_n F_n(x)$ in $C(\mathbb{T})$.
Further, since every $F_n(x)$ is Lipschitz-continuous with Lipschitz-constant $K = \sup_n \|f_n\|_{\infty}$, so is their uniform limit $F_{\eta}(x)$. Thus $F_{\eta}(x)$ is differentiable a.e on $\mathbb{T}$. Then we have the following result.

\begin{lem}
$F_{\eta}(x) = G_{\eta}(x)$ holds for every $\eta \in \mathbb{N}_0^*$. In other words, $F^{\prime}_{\eta}(x) = f_{\eta}(x)$ m-a.e on $\mathbb{T}$.
\end{lem}

\begin{prf}
Since the definition of weak*-convergence, for any $x \in [0, 1]$ we have 
\[F_{\eta}(x) = \eta\mathchar`-\lim_n F_n(x) = \eta\mathchar`-\lim_n \int_0^x f_n(t)dt = \int_0^x f_{\eta}(t)dt = G_{\eta}(x). \] 
Hence $F^{\prime}_{\eta}(x) = f_{\eta}(x)$ also holds by the Lebesgue differentiation theorem.
\end{prf}

\begin{thm}
For each $f \in L^{\infty}(\mathbb{R}_+)$ it holds that
\[\overline{K}(f) = \sup_{\eta \in \mathbb{N}^*} \esssup_{x \in [0, 1]} \tilde{f}_{\eta}(x). \]
\end{thm}

\begin{prf}
Let $f(x) \in L^{\infty}(\mathbb{R})$ and let us define
\[\tilde{F}_n(x) = \int_0^x \tilde{f}_n(t)dt, \quad x \in [0, 1], \]
\[\tilde{F}_{\eta}(x) = \eta\mathchar`-\lim_n \tilde{F}_n(x), \quad x \in [0, 1]. \]
Remark that by the Lemma 3.1, it holds that
\[\tilde{F}_{\eta}(x) = \int_0^x \tilde{f}_{\eta}(t)dt, \ i.e., \ \tilde{F}_{\eta}^{\prime}(x) = \tilde{f}_{\eta}(x), \quad x \in [0, 1]. \]
Then we have 
\begin{align}
\overline{K}(f) &= \lim_{\theta \to 0^+} \limsup_{x \to \infty} \frac{1}{\theta} \int_x^{x+\theta} f(t)dt \notag \\
&= \lim_{\theta \to 0^+} \limsup_{n \to \infty} \sup_{x \in [0, 1-\theta]} \frac{1}{\theta} \int_x^{x+\theta} \tilde{f}_n(t)dt \notag \\
&= \lim_{\theta \to 0^+} \limsup_{n \to \infty} \sup_{x \in [0, 1-\theta]} \frac{\tilde{F}_n(x+\theta) - \tilde{F}_n(x)}{\theta} \notag \\
&= \lim_{\theta \to 0^+} \sup_{\eta \in \mathbb{N}^*} \sup_{x \in [0, 1-\theta]} \eta\mathchar`-\lim_n \frac{\tilde{F}_n(x+\theta) - \tilde{F}_n(x)}{\theta} \notag \\
&= \lim_{\theta \to 0^+} \sup_{\eta \in \mathbb{N}^*} \sup_{x \in [0, 1-\theta]} \frac{\tilde{F}_{\eta}(x+\theta) - \tilde{F}_{\eta}(x)}{\theta} \notag \\
&\ge \sup_{\eta \in \mathbb{N}^*} \sup_{x \in [0, 1)} \limsup_{\theta \to 0^+}  \frac{\tilde{F}_{\eta}(x+\theta) - \tilde{F}_{\eta}(x)}{\theta} \notag \\
&\ge \sup_{\eta \in \mathbb{N}^*} \esssup_{x \in [0, 1)} \limsup_{\theta \to 0^+}  \frac{\tilde{F}_{\eta}(x+\theta) - \tilde{F}_{\eta}(x)}{\theta} \notag \\
&= \sup_{\eta \in \mathbb{N}^*} \esssup_{x \in [0, 1]} \tilde{f}_{\eta}(x). \notag
\end{align}
On the other hand, for any $\eta \in \mathbb{N}^*$, $\theta > 0$ and $x \in [0, 1-\theta]$, we have
\begin{align}
\frac{\tilde{F}_{\eta}(x + \theta) - \tilde{F}_{\eta}(x)}{\theta} &= \frac{1}{\theta} \int_x^{x + \theta} \tilde{f}_{\eta}(t)dt
\le \frac{1}{\theta} \cdot \theta \cdot \esssup_{x \in [0, 1]} \tilde{f}_{\eta}(x) \notag \\
&= \esssup_{x \in [0, 1]} \tilde{f}_{\eta}(x). \notag
\end{align}
Thus we have
\[\lim_{\theta \to 0^+} \sup_{\eta \in \mathbb{N}^*} \sup_{x \in [0, 1-\theta]} \frac{\tilde{F}_{\eta}(x + \theta) - \tilde{F}_{\eta}(x)}{\theta} \le \sup_{\eta \in \mathbb{N}^*} \esssup_{x \in [0, 1]} \tilde{f}_{\eta}(x). \] 
\[\therefore \overline{K}(f) = \sup_{\eta \in \mathbb{N}^*} \esssup_{x \in [0, 1]} \tilde{f}_{\eta}(x). \]
\end{prf}

Now we define for each $\omega = (\eta, t) \in \Omega^*$ a linear operator $T_{\omega} : L^{\infty}(\mathbb{R}_+) \rightarrow L^{\infty}(\mathbb{R})$ as follows; let us consider a set of functions $\{f(x+s)\}_{s \ge 0}$ in $L^{\infty}(\mathbb{R}_+)$. This is bounded and weak* relatively compact set of $L^{\infty}(\mathbb{R}_+)$. Thus for each $\omega \in \Omega^*$ we can define its limit along $\omega$, i.e., $\omega\mathchar`-\lim_s f(x+s)$ in $L^{\infty}(\mathbb{R}_+)$ with respect to the weak* topology.  Moreover, we can extend it to the negative direction by
\[(T_{\omega}f)(x) = f_{\tau^{-N}\omega}(N+x), \quad x \in [-N, 0], \]
for every $N > 0$. Then it is easy to see that this definition of $(T_{\omega}f)(x)$ is equal to the following.
\[(T_{\omega}f)(x) = \tilde{f}_{\tau^{[x+t]}\eta}(x+t-[x+t]), \quad x \in \mathbb{R}. \]
Notice that for $f(x) \in C_{ub}(\mathbb{R}_+)$, it holds that $(Tf)_{\omega}(x) = \overline{f}(\tau^x\omega)$, i.e., the restriction of $\overline{f}(\omega)$ to the orbit $o(\omega)$ of $\omega$.

Thus Theorem 3.1 can be expressed as follows.
\begin{cor}
$\overline{K}(f) = \sup_{\omega \in \Omega^*} \esssup_{x \in \mathbb{R}} (T_{\omega}f)(x)$.
\end{cor}

For a fucntion $f \in L^{\infty}(\mathbb{R}_+)$, let us denote $K(f) = \alpha$ if $\varphi(f) = \alpha$ for every $\varphi \in \mathcal{K}$. Then the following corollary follows immediately by the above corollary.
\begin{cor}
$K(f) = \alpha$ if and only if $w^*\mathchar`-\lim_s f(x+s) = \alpha$, where the symbol $w^*\mathchar`-\lim$ represents weak*-convergence in $L^{\infty}(\mathbb{R}_+)$.
\end{cor}

Now the following result is important for our intention. For the sake of simplicity, for any $f \in L^{\infty}(\mathbb{R}_+)$, let us denote $(T_{\omega}f)(x) = f_{\omega}(x)$ and 
\[(Sf)(x) = \frac{1}{e^x} \int_0^x f(t)e^tdt, \quad x \ge 0. \]

\begin{thm}
$R(f) = \alpha$ if and only if $K(f) = \alpha$.
\end{thm}

\begin{prf}
First we observe that
\begin{align}
R(f) = \alpha &\Longleftrightarrow \lim_{x \to \infty} (Sf)(x) = \alpha  \notag  \\
&\Longleftrightarrow (\overline{Sf})(\omega) = \alpha \ on \ \Omega^* \notag \\
&\Longleftrightarrow (Sf)_{\omega}(x) = \alpha \ for \ every \ \omega \in \Omega^*. \notag
\end{align}

For each $f(x) \in L^{\infty}(\mathbb{R}_+)$, notice that
\[f(x) = (Sf)(x) + (Sf)^{\prime}(x), \quad x \ge 0. \]
Operating $T_{\omega}$ both sides and recalling that $((Sf)^{\prime})_{\omega}(x) = ((Sf)_{\omega})^{\prime}(x)$ by Lemma 3.1, we have 
\[f_{\omega}(x) = (Sf)_{\omega}(x) +((Sf)_{\omega})^{\prime}(x), \quad x \in \mathbb{R}. \]
Hence
\[(e^x \cdot (Sf)_{\omega}(x))^{\prime} = e^x \cdot f_{\omega}(x), \quad x \in \mathbb{R}. \]
Since $(Sf)_{\omega}(x)$ is bounded on $\mathbb{R}$, we obtain immediately that for each $\omega \in \Omega^*$, $(Sf)_{\omega}(x) = \alpha$ if and only if $f_{\omega}(x) = \alpha$. This completes the proof.
\end{prf}

\begin{cor}
$R(f) = \alpha$ if and only if $w^*\mathchar`-\lim_s f(x+s) = \alpha$.
\end{cor}

Hnece we have $F_0 = \{f(x) \in L^{\infty}(\mathbb{R}_+) : w^*\mathchar`-\lim_s f(x+s) = 0\}$.
Then by Corollary 3.1 it is immediate that $\mathcal{K} \subseteq \mathcal{D}$. In fact, we can show the reverse inclusion:

\begin{thm}
$\mathcal{K} = \mathcal{D}$ holds.
\end{thm}

\begin{prf}
It is sufficient to show that $\mathcal{D} \subseteq \mathcal{K}$ and by the Krein-Milman theorem it is equivalent to the assertion that $\overline{K}_0(f) \le \overline{K}(f)$ for every $f \in L^{\infty}(\mathbb{R}_+)$, where
\[\overline{K}_0(f) = \sup_{\varphi \in S^+_{L^{\infty}(\mathbb{R}_+)^*} \cap F_0^{\bot}} \varphi(f). \]

 First, remark that
\begin{align}
\overline{K}_0(f) &= \inf_{h \in F_0} \esssup_{x \in [0, \infty)} (f(x) - h(x))   \notag \\
&= \inf_{h \in F_0} \sup_{n \in \mathbb{N}^0} \esssup_{x \in [0, 1]} (\tilde{f}_n(x) -\tilde{h}_n(x)). \notag
\end{align}                 

Now for each $h \in F_0$ we put
\[H_n(x) = \int_0^x (\tilde{f}_n(t) -\tilde{h}_n(t))dt, \quad x \in [0, 1], \ n=0,1,2, \ldots. \]

Then by the assumption that $h(x) \in F_0$ we have
\[H_n(x) \rightarrow F_n(x) \quad uniformly \ in \ x \in [0, 1] \ as \ n \to \infty. \]

Also we note
\[\tilde{f}_n(x) -\tilde{h}_n(x) = \lim_{\theta \to 0} \frac{H_n(x+\theta) - H_n(x)}{\theta} \quad m\mathchar`-a.e. \ on \ [0, 1], \ 
n=0,1,2, \ldots. \]

Hence we get that
\begin{align}
\overline{K}_0(f) &= \inf_{h \in F_0} \sup_{n \in \mathbb{N}^0} \sup_{x \in [0, 1]} \limsup_{\theta \to 0} \frac{H_n(x+\theta) - H_n(x)}{\theta} \notag \\
&\le \inf_{h \in F_0} \limsup_{\theta \to 0} \sup_{n \in \mathbb{N}^0} \sup_{x \in [0, 1-\theta]} \frac{H_n(x+\theta) - H_n(x)}{\theta} \notag \\
&= \inf_{h \in F_0} \limsup_{\theta \to 0} \limsup_{n \in \mathbb{N}^0} \sup_{x \in [0, 1-\theta]} \frac{H_n(x+\theta) - H_n(x)}{\theta} \notag \\
&= \inf_{h \in F_0} \limsup_{\theta \to 0} \limsup_{n \in \mathbb{N}^0} \sup_{x \in [0, 1-\theta]} \frac{F_n(x+\theta) - F_n(x)}{\theta} \notag \\
&= \limsup_{\theta \to 0} \limsup_{n \in \mathbb{N}^0} \sup_{x \in [0, 1-\theta]} \frac{1}{\theta} \int_{n+x}^{n+x+\theta} f(t)dt \notag \\
&= \limsup_{\theta \to 0} \limsup_{x \to \infty} \frac{1}{\theta} \int_x^{x+\theta} f(t)dt = \overline{K}(f) \notag 
\end{align}
\end{prf}

\begin{thm}
$\sup_{\varphi \in \mathcal{D}} \varphi(f) = \overline{K}(f)$ holds for each $f \in L^{\infty}(\mathbb{R}_+)$.
\end{thm}

\section{Extremal values of $\mathcal{C}$}
We consider the set of continuous linear functionals $\psi$ on $L^{\infty}(\mathbb{R}_+^{\times})$ for which 
\[\psi(f) \le \overline{P}(f) \]
holds for every $f \in L^{\infty}(\mathbb{R}_+^{\times})$ and denote it by $\mathcal{P}$.
In this section we deal with the multiplicative version of the preceding section. First, for given $f(x) \in L^{\infty}(\mathbb{R}_+^{\times})$, we define in a similar way as the former section the sequence $\{\hat{f}_n(x)\}_{n \ge 0}$ of elements of $L^{\infty}([1, e])$ by
\[\hat{f}_n(x) = f(e^n x), \quad x \in [1, e], \ n = 0, 1, 2, \cdots.   \]
Similarly we consider its limit along $\eta \in \mathbb{\beta}\mathbb{N}^0$, denoted by $\hat{f}_{\eta}(x) = \eta\mathchar`-\lim_n \hat{f}_n(x)$, in $L^{\infty}([1, e])$ with respect to its weak*-topology. Also let us define
\[\hat{F}_n(x) = \int_1^x \hat{f}_n(t)dt \in C([1, e]), \quad n=0,1,2, \cdots . \]
\[\hat{F}_{\eta}(x) = \eta\mathchar`-\lim_n \hat{F}_n(x), \quad \eta \in \mathbb{N}^*. \]
Then we have an analogous result to Lemma 3.1:

\begin{lem}
$\hat{F}^{\prime}_{\eta}(x) = \hat{f}_{\eta}(x)$ m-a.e. on $[1, e]$.
\end{lem}

Then we can show an analogue of Theorem 3.1 by similar arguments:

\begin{thm} For every $f \in L^{\infty}(\mathbb{R}_+^{\times})$ it holds that
\[\overline{P}(f) = \sup_{\eta \in \mathbb{N}^*} \esssup_{x \in [1, e]} \hat{f}_{\eta}(x). \]
\end{thm}

\begin{prf}
For given $f(x) \in L^{\infty}(\mathbb{R}_+^{\times})$, we have

\begin{align}
\overline{P}(f) &= \lim_{\theta \to 1^+} \limsup_{x \to \infty} \frac{1}{\theta x - x} \int_x^{\theta x} f(t)dt \notag \\
&= \lim_{\theta \to 1^+} \limsup_{n \to \infty} \sup_{x \in [1, e/\theta]} \frac{1}{\theta x-1} \int_x^{\theta x} \hat{f}_n(t)dt    \notag \\
&= \lim_{\theta \to 1^+} \limsup_{n \to \infty} \sup_{x \in [1, e/\theta]} \frac{\hat{F}_n(\theta x) -\hat{F}_n(x)}{\theta x - x}    \notag \\
&= \lim_{\theta \to 1^+} \sup_{\eta \in \mathbb{N}^*} \sup_{x \in [1, e/\theta]} \eta\mathchar`-\lim_n \frac{\hat{F}_n(\theta x) -\hat{F}_n(x)}{\theta x - x}   \notag \\
&= \lim_{\theta \to 1^+} \sup_{\eta \in \mathbb{N}^*} \sup_{x \in [1, e/\theta]} \frac{\hat{F}_{\eta}(\theta x) - \hat{F}_{\eta}(x)}
{\theta x -x}   \notag \\
&\ge \sup_{\eta \in \mathbb{N}^*} \sup_{x \in [1,e]} \limsup_{\theta \to 1^+} \frac{\hat{F}_{\eta}(\theta x) - \hat{F}_{\eta}(x)}
{\theta x -x}   \notag 
\end{align}
\begin{align}
&\ge \sup_{\eta \in \mathbb{N}^*} \esssup_{x \in [1,e]} \limsup_{\theta \to 1^+} \frac{\hat{F}_{\eta}(\theta x) - \hat{F}_{\eta}(x)}
{\theta x -x}   \notag \\
&= \sup_{\eta \in \mathbb{N}^*} \esssup_{x \in [1,e]} \hat{f}_{\eta}(x). \notag
\end{align}
On the other hand, for any $\eta \in \mathbb{N}^*, \theta > 1$ and $x \in [1, e/\theta]$, we have 
\begin{align}
\frac{\hat{F}_{\eta}(\theta x) - \hat{F}_{\eta}(x)}{\theta x -x} = \frac{1}{\theta x - x} \int_x^{\theta x} \hat{f}_{\eta}(t)dt
&\le \frac{1}{\theta x - x} \cdot (\theta x - x) \cdot \esssup_{x \in [1,e]} \hat{f}_{\eta}(x) \notag \\
&= \esssup_{x \in [1,e]} \hat{f}_{\eta}(x). \notag
\end{align}
Thus we get
\[\lim_{\theta \to 1^+} \sup_{\eta \in \mathbb{N}^*} \sup_{x \in [1, e/\theta]} \frac{\hat{F}_{\eta}(\theta x) - \hat{F}_{\eta}(x)}
{\theta x -x} \le \sup_{\eta \in \mathbb{N}^*} \esssup_{x \in [1,e]} \hat{f}_{\eta}(x). \]
\[\therefore \overline{P}(f) = \sup_{\eta \in \mathbb{N}^*} \esssup_{x \in [1,e]} \hat{f}_{\eta}(x).  \]
\end{prf}

\bigskip

Similarly as $\overline{K}(f)$, we define for each $\omega \in \Omega$ a linear operator $P_{\omega} : L^{\infty}(\mathbb{R}^{\times}) \longrightarrow L^{\infty}(\mathbb{R}^{\times})$ as follows; let us consider the set of functions $\{f(rx)\}_{r \ge 1}$ in $L^{\infty}(\mathbb{R}_+^{\times})$. This is bounded and weak* relatively compact set of $L^{\infty}(\mathbb{R}_+^{\times})$. Thus for each $\omega \in \Omega^*$ we can define its limit along the ultrafilter $e^\omega = \{e^A : A \in \omega \}$, i.e., $e^{\omega}\mathchar`-\lim_r f(rx) = \omega\mathchar`-\lim_r f(e^rx)$. Remark that in this way one obtains just a function in $L^{\infty}(\mathbb{R}_+^{\times})$, but can extend it to the whole space $\mathbb{R}^{\times}$ in the same way as the definiton of $T_{\omega}$. Then we have 
\[(P_{\omega}f)(x) = \hat{f}_{\tau^{[\log x\theta]}\eta}(x\theta/e^{[\log x\theta]}), \quad x \in \mathbb{R}^{\times}, \]
where $\omega = (\eta, t)$ and $\theta = e^t$. Hence Theorem 3.1 can be expressed as follows.

\begin{cor}
$\overline{P}(f) = \sup_{\omega \in \Omega^*} \esssup_{x \in \mathbb{R}} (P_{\omega}f)(x)$.
\end{cor}

Now we take up the relation between $\tilde{f}_{\eta}(x) \in L^{\infty}([0, 1])$ defined in the former section and $\hat{f}_{\eta}(x) \in L^{\infty}([1, e])$. We define a linear operator $W_0$ by
\[W_0 : L^{\infty}([1, e]) \longrightarrow L^{\infty}([0, 1]), \quad (W_0f)(x) = f(e^x). \]
Then notice that
\[(W_0\hat{f}_n)(x) = (\widetilde{W_0f})_n(x), \quad n=0,1,2, \ldots, \]
holds for each $f \in L^{\infty}(\mathbb{R}_+^{\times})$. Further, we have the following result.

\begin{lem}
$(W_0\hat{f}_{\eta})(x) = (\widetilde{W_0f})_{\eta}(x)$ for each $\eta \in \mathbb{N}^*$ and $f \in L^{\infty}(\mathbb{R}_+^{\times})$.
\end{lem}

\begin{prf}
For any $f \in L^{\infty}(\mathbb{R}_+^{\times})$ and $\eta \in \mathbb{N}^*$, by the definition of the weak* topology,
\[
\eta\mathchar`-\lim_n \hat{f}_n(x) = \hat{f}_{\eta}(x)  
\Longleftrightarrow \eta\mathchar`-\lim_n \int_1^e \hat{f}_n(t)\phi(t)dt = \int_1^e \hat{f}_{\eta}(t)\phi(t)dt 
\]
for every $\phi \in L^1([1,e])$. Notice that by integration by substitution for any $g(x) \in L^1([1, e])$ we have
\[\int_1^e g(t)dt = \int_0^1 g(e^t)\cdot e^tdx. \]
Hence for every $\phi \in L^1([1,e])$, we have
\[\eta\mathchar`-\lim_n \int_0^1 \hat{f}_n(e^t)\phi(e^t)e^tdt = \int_0^1 \hat{f}_{\eta}(e^t)\phi(e^t)e^tdt. \]
Namely, it holds that
\[\eta\mathchar`-\lim_n \int_0^1 (\widetilde{W_0f})_n(t)\psi(t)dt = \int_0^1 (W_0\hat{f}_{\eta})(t)\psi(t)dt \]
where $\psi(x) = \phi(e^x)e^x$.
Notice that for every $\psi(x) \in L^1([0,1])$, $\phi(y) = \frac{\psi(\log y)}{y} \in L^1([1,e])$ and $\phi(e^x)e^x = \psi(x)$ holds. This means that
\[\eta\mathchar`-\lim_n \int_0^1 (\widetilde{W_0f})_n(t)\psi(t)dt = \int_0^1 (W_0\hat{f}_{\eta})(t)\psi(t)dt \]
for every $\psi(x) \in L^1([0,1])$.
Therefore, $(\widetilde{W_0f})_{\eta}(x) = \eta\mathchar`-\lim_n (\widetilde{W_0f})_n(x) = (W_0\hat{f}_{\eta})(x)$ in $L^{\infty}([0,1])$ for every $\eta \in \mathbb{N}^*$. We complete the proof.
\end{prf}

It is possible to extend this result to the relation between $(T_{\omega}f)(x)$ and $(P_{\omega}f)(x)$ as follows; let us define a linear operator $W$ by
\[W : L^{\infty}(\mathbb{R}_+^{\times}) \longrightarrow L^{\infty}(\mathbb{R}), \quad (Wf)(x) = f(e^x). \]

\begin{cor}
For every $f \in L^{\infty}(\mathbb{R}_+^{\times})$ and $\omega \in \Omega$, $T_{\omega}Wf = WP_{\omega}f$ holds. Namely, $P_{\omega} = W^{-1}T_{\omega}W$ holds.
\end{cor}

\begin{prf}
Let $x \in \mathbb{R}$, $\omega = (\eta, t) \in \Omega^*$, $\theta = e^t$ and $f \in L^{\infty}(\mathbb{R}_+^{\times})$. Then we have
\begin{align}
(WP_{\omega}f)(x) &= \hat{f}_{\tau^{[\log{e^x} \cdot \theta]}\eta} (e^x \cdot \theta / e^{[\log e^x \cdot \theta]}) \notag \\
&= \hat{f}_{\tau^{[x + \log\theta]}\eta} (e^x \cdot \theta / e^{[x + \log\theta]}) \notag \\
&= \hat{f}_{\tau^{[x+t]}\eta} (e^{x+t} / e^{[x+t]}) \notag \\
&= \hat{f}_{\tau^{[x+t]}\eta} (e^{x+t - [x+t]})  \notag \\
&= (W\hat{f}_{\tau^{[x+t]}\eta}) (x+t-[x+t]) \notag \\
&= (\widetilde{Wf})_{\tau^{[x+t]}\eta} (x+t-[x+t]) \notag \\
&= (T_{\omega}Wf)(x). \notag
\end{align} 
\end{prf}

\medskip

Remark that for each $\omega \in \Omega^*$, $\esssup_{x \in \mathbb{R}^{\times}} (P_{\omega}f)(x) = \esssup_{x \in \mathbb{R}} (P_{\omega}f)(e^x) = \esssup_{x \in \mathbb{R}} (T_{\omega}Wf)(x)$ holds and then we have the following result.

\begin{thm}
$\overline{P}(f) = \overline{K}(Wf)$ for every $f \in L^{\infty}(\mathbb{R}_+^{\times})$.
\end{thm}

Let us consider the adjoint operator $W^*$ of $W$, which is a linear isomeory from $L^{\infty}(\mathbb{R}_+)^*$ onto $L^{\infty}(\mathbb{R}_+^{\times})^*$. Theorem 4.2 shows that $\mathcal{K}$ and $\mathcal{P}$ are affinely homeomorphic via $W^*$. 

For a function $f \in L^{\infty}(\mathbb{R}_+^{\times})$, we denote $P(f) = \alpha$ if $\psi(f) = \alpha$ for every $\psi \in \mathcal{P}$.

\begin{thm}
Let $f \in L^{\infty}(\mathbb{R}_+^{\times})$. Then $M(f) = \alpha$ if and only if $P(f) = \alpha$.
\end{thm}

\begin{prf}
By Theorem 4.2, $P(f) = \alpha$ is equivalent to $K(Wf) = \alpha$. Also, notice that
\[\lim_{x \to \infty} \frac{1}{x} \int_1^x f(t)dt = \lim_{x \to \infty} \frac{1}{x} \int_0^{\log x} (Wf)(t)e^{t}dt = \lim_{x \to \infty} \frac{1}{e^x} \int_0^x (Wf)(t)e^tdt,  \]
which shows that $M(f) = \alpha$ is equivalent to $R(Wf) = \alpha$. Thus, by Theorem 3.2, we have
\[M(f) = \alpha \Longleftrightarrow R(Wf) = \alpha \Longleftrightarrow K(Wf) = \alpha \Longleftrightarrow P(f) = \alpha. \]
\end{prf}

\begin{cor}
$M(f) = \alpha$ if and only if $w^*\mathchar`-\lim_r f(rx) = \alpha$.
\end{cor}

As we have seen in the above proof, $f(x) \in E_0$ is equivalent to $(Wf)(x) \in F_0$, which menas that $WE_0 = F_0$, and thus $W^*\mathcal{D} = \mathcal{C}$. Therefore since we have already shown that $W^*\mathcal{K} = \mathcal{P}$ and $\mathcal{K} = \mathcal{D}$, we get $\mathcal{P} = \mathcal{C}$. Now the following result of the extremal value attined by $\mathcal{C}$ follows immediately.

\begin{thm}
$\sup_{\psi \in \mathcal{C}} \psi(f) = \overline{P}(f)$ for each $f \in L^{\infty}(\mathbb{R}_+^{\times})$.
\end{thm}

Finally, we will give another expression of $\sup_{\psi \in \mathcal{C}} \psi(f)$. We define a sublinear functional $\overline{Q}$ on $L^{\infty}(\mathbb{R}_+^{\times})$ as follows.
\[\overline{Q}(f) = \lim_{\theta \to 1^+} \limsup_{x \to \infty} \frac{1}{\log \theta} \int_x^{\theta x} f(t) \frac{dt}{t}, \]
where $f(x) \in L^{\infty}(\mathbb{R}_+^{\times})$. 

\begin{thm}
$\overline{Q}(f) = \overline{K}(Wf)$ holds for every $f \in L^{\infty}(\mathbb{R}_+^{\times})$.
\end{thm}

\begin{prf}
For any $f \in L^{\infty}(\mathbb{R}_+^{\times}), x \ge 0$ and $\theta > 0$, we have that
\begin{align}
\frac{1}{\theta} \int_x^{x + \theta} (Wf)(t)dt &= \frac{1}{\theta} \int_x^{x+\theta} f(e^t)dt \notag \\
&= \frac{1}{\theta} \int_{e^x}^{e^x \cdot e^{\theta}} f(s) \frac{1}{ds}  \notag 
\end{align}
Now we put $r=e^x$, then $r$ tends to $\infty$ as $x$ tends to $\infty$. And then put $y=e^{\theta}$ and then $y$ tends to
$1^+$ as $\theta$ tends to $0^+$. Hence we have 
\[\overline{K}(Wf) = \lim_{y \to 1^+} \limsup_{r \to \infty} \frac{1}{\log y} \int_r^{ry} f(s) \frac{ds}{s} = \overline{Q}(f). \]
\end{prf}

Hence by Theorem 4.2 and Theorem 4.4, we get the following result.
\begin{cor}
$\sup_{\psi \in \mathcal{C}} \psi(f) = \overline{Q}(f)$ for each $f \in L^{\infty}(\mathbb{R}_+^{\times})$.
\end{cor}

\medskip

\section{Applications to Ces\`{a}ro mean}
In this section, we deal with the discrete version of the preceding results, that is, normalized positive linear functionals on $l_{\infty}$ which extend Ces\`{a}ro mean. First, we consider the linear operators $V$ and $V_1$ defined as follows.
\[V : l_{\infty} \longrightarrow L^{\infty}(\mathbb{R}_+^{\times}), \quad (Vf)(x) = f([x]). \]
and
\[V_1 : L^{\infty}(\mathbb{R}_+^{\times}) \longrightarrow l_{\infty}, \quad (V_1f)(n) = \int_n^{n+1} f(t)dt. \]
Let $V^*$ and $V_1^*$ be their adjoint operators respectively. Then we have the following result.

\begin{thm}
$\mathcal{C}$ and $\mathcal{C}_d$ are affinely homeomorphic via $V^*$.
\end{thm}

\begin{prf}
First, notice that for each $f \in l_{\infty}$, we have
\[\frac{1}{n} \sum_{i=1}^n f(i) = \frac{1}{n} \int_1^{n+1} (Vf)(t)dt. \]
Hence if $\varphi \in \mathcal{C}$ and $f \in l_{\infty}$ with $M_d(f)$ exists, then
\[(V^*\varphi)(f) = \varphi(Vf) = M(Vf) = M_d(f). \]
Thus $V^*\varphi \in \mathcal{C}_d$ holds for every $\varphi \in \mathcal{C}$. Next we show that $V^*$ is one to one. Suppose that $\varphi_1, \varphi_2 \in \mathcal{C}$ with $\varphi_1 \not= \varphi_2$. Then there is some $f \in L^{\infty}(\mathbb{R}_+^{\times})$ such that $\varphi_1(f) \not= \varphi_2(f)$. On the other hand, for every $\varphi \in \mathcal{C}$,
it holds that $\varphi(f) = \varphi(VV_1f)$ since
\[\lim_{x \to \infty} \frac{1}{x} \int_1^x (f(t) - (VV_1f)(t))dt = 0. \]
Then we have 
\[\varphi_1(f) = \varphi_1(VV_1f) = (V^*\varphi_1)(V_1f), \quad \varphi_2(f) = \varphi_2(VV_1f) = (V^*\varphi_2)(V_1f) \]
Thus, by the assumption that $\varphi_1(f) \not= \varphi_2(f)$, we have $(V^*\varphi_1)(V_1f) \not= (V^*\varphi_2)(V_1f)$.
Hence $V^*\varphi_1 \not= V^*\varphi_2$. This shows that $V^* : \mathcal{C} \rightarrow \mathcal{C}_d$ is one to one.

Next we show that $V^* : \mathcal{C} \rightarrow \mathcal{C}_d$ is onto. Given any $\psi \in \mathcal{C}_d$, let us take $\varphi = V_1^*\psi \in \mathcal{C}$. Remark that $V_1Vf = f$ for every $f \in l_{\infty}$ and then we have $V^*\varphi = V^*V_1^*\psi = \psi$. Hence $V^* : \mathcal{C} \rightarrow \mathcal{C}_d$ is onto. We obtain the result.
\end{prf}

\begin{lem}
For each $f \in l_{\infty}$, it holds that
\[\lim_{\theta \to 1^+} \limsup_{n \to \infty} \frac{1}{\theta n - n} \int_n^{\theta n} (Vf)(t)dt = \lim_{\theta \to 1^+} \limsup_{n \to \infty} \frac{1}{\theta n - n} \sum_{i \in [n, \theta n]} f(i). \]
\end{lem}

\begin{prf}
Let us $n \in \mathbb{N}$, $\theta > 1$ and $f \in l_{\infty}$. The assertion follows from the following equation.
\[\int_n^{[\theta n]+1} (Vf)(t)dt = \sum_{i=n}^{[\theta n]} f(i). \]
\end{prf}

\begin{thm}
$\sup_{\psi \in \mathcal{C}_d} \psi(f) = \overline{P}_d(f)$ holds for each $f \in l_{\infty}$.
\end{thm}

\begin{prf}
We have that
\begin{align}
\sup_{\psi \in \mathcal{C}_d} \psi(f) &= \sup_{\varphi \in \mathcal{C}} (V^*\varphi)(f) = \sup_{\varphi \in \mathcal{C}} \varphi(Vf) = \overline{P}(Vf) \notag \\
&= \lim_{\theta \to 1^+} \limsup_{x \to \infty} \frac{1}{\theta x - x} \int_x^{\theta x} (Vf)(t)dt \notag \\
&= \lim_{\theta \to 1^+} \limsup_{n \to \infty} \frac{1}{\theta n - n} \int_n^{\theta n} (Vf)(t)dt \notag \\
&= \lim_{\theta \to 1^+} \limsup_{n \to \infty} \frac{1}{\theta n - n} \sum_{i \in [n, \theta n]} f(i). \notag \\
&= \overline{P}_d(f). \notag
\end{align}
\end{prf}

We define a sublinear functional $\overline{Q}_d$ on $l_{\infty}$, which is a discrete version of the sublinear functional $\overline{Q}$ as follows:
\[Q_d(f) =  \lim_{\theta \to 1^+} \limsup_{n \to \infty} \frac{1}{\log \theta} \sum_{i \in [n, \theta n]} \frac{f(i)}{i}, \quad f \in l_{\infty}. \]
\begin{lem}
For each $f \in l_{\infty}$, it holds that
\[\lim_{\theta \to 1^+} \limsup_{n \to \infty} \frac{1}{\log \theta} \int_n^{\theta n} (Vf)(t) \frac{dt}{t} = \lim_{\theta \to 1^+} \limsup_{n \to \infty} \frac{1}{\log \theta} \sum_{i \in [n, \theta n]} \frac{f(i)}{i}. \]
\end{lem}

\begin{prf}
Let us $n \in \mathbb{N}$, $\theta > 1$ and $f \in l_{\infty}$. The lemma follows immediately from the following computation.
\begin{align}
&\Bigg|\int_n^{[\theta n] + 1} (Vf)(t) \frac{dt}{t} - \sum_{i \in [n, \theta n]} \frac{f(i)}{i} \Bigg| \notag \\
&= \Bigg|f(n) \cdot \left(\frac{1}{n} - \log\left(1+\frac{1}{n}\right)\right) + f(n+1) \cdot \left(\frac{1}{n+1} - \log\left(1+\frac{1}{n+1}\right)\right) + \ldots \notag \\
&+ f([\theta n]) \cdot \left(\frac{1}{[\theta n]} - \log\left(1+\frac{1}{[\theta n]}\right)\right)\Bigg| \notag \\
&\le |f(n)| \cdot \frac{1}{2} \cdot \frac{1}{n^2} + |f(n+1)| \cdot \frac{1}{2} \cdot \frac{1}{(n+1)^2} + \ldots |f([\theta n])| \cdot \frac{1}{2} \cdot \frac{1}{[\theta n]^2} \notag \\
&\le \|f\|_{\infty} \cdot \frac{\pi^2}{12}. \notag
\end{align}
\end{prf}

Then we can show the following result in a similar way as Theorem 5.2.
\begin{thm}
$\sup_{\psi \in \mathcal{C}_d} \psi(f) = \overline{Q}_d(f)$ holds for each $f \in l_{\infty}$.
\end{thm}

\end{document}